\documentclass[11pt, reqno]{amsart}

\usepackage{fullpage}   
\usepackage{amssymb}    
\usepackage{url}

\theoremstyle{plain}

\newtheorem{theorem}{Theorem}
\newtheorem{proposition}{Proposition}
\newtheorem{prob}{Problem}

\theoremstyle{remark}
\newtheorem*{acknowledgment}{Acknowledgment}

\numberwithin{equation}{section}

\newcommand{\seclabel}[1]{\label{sec:#1}}   

\title{A 2-Base for Inverse Semigroups}

\author{Jo\~{a}o Ara\'{u}jo}
\author{Michael Kinyon}
\author{R. Padmanabhan}

\address[Ara\'{u}jo]{Centro de \'{A}lgebra \\
Universidade de Lisboa \\
1649-003 Lisboa \\ Portugal\\
and\\
Universidade Aberta\\
1269--001 Lisboa \\ Portugal}
\email{\url{jaraujo@ptmat.fc.ul.pt}}

\address[Kinyon]{Department of Mathematics \\
University of Denver \\ 2360 S Gaylord St \\ Denver, Colorado 80208 USA}
\email{\url{mkinyon@du.edu}}

\address[Padmanabhan]{Department of Mathematics \\
University of Manitoba\\ Winnipeg R3T 2N2 \\ CANADA}
\email{\url{padman@cc.umanitoba.ca}}

\begin{document}
\maketitle
\vspace{-0.5cm}
\begin{center}
{\em Dedicated to the memory of Bill McCune (1953-2011)}
\end{center}

\begin{abstract}
    An open problem in the theory of inverse semigroups was whether the variety of such semigroups, when
viewed as algebras with a binary operation and a unary operation, is $2$-based, that is, has a base for
its identities consisting of $2$ independent axioms. In this note, we announce the affirmative
solution to this problem: the identities
\[
\quad x(x'x) = x \qquad
\quad x  (x'  (y  (y'  ((z  u)'  w')'))) = y  (y'  (x  (x'  ((w  z)  u))))
\]
form a base for inverse semigroups where ${}'$ turns out to be the natural inverse operation. We recount here the history of the problem including our previous efforts to find a $2$-base using automated
edduction and the method that finally worked. We describe our efforts to simplify the proof using \textsc{Prover9}, present the simplified proof itself and conclude with some open problems. A humanized proof will be submitted elsewhere.
\end{abstract}

\section{Introduction}
\seclabel{intro}

The notion of inverse in semigroup theory generalizes the corresponding notion in group theory.
Given an element $a$ of a semigroup, an element $b$ is said to be an \emph{inverse} of $a$ if the
equations $aba = a$ and $bab = b$ hold. A semigroup in which every element has an inverse is said to
be \emph{regular}. A regular semigroup in which each element has a \emph{unique} inverse is called an \emph{inverse semigroup}. Arguably, groups, inverse semigroups and regular semigroups constitute the most important classes of semigroups. A standard reference for semigroup theory in general is \cite{Howie}, and another for inverse semigroups in particular is \cite{Petrich}.

As idempotents shape the structure of regular semigroups to a large extent, it is no surprise that these classes can be defined in terms of idempotents, that is,
\begin{itemize}
\item  a group is a regular semigroup with exactly one idempotent;
\item  an inverse semigroup is a regular semigroup in which the idempotents commute.
\end{itemize}

Unlike general regular semigroups, inverse semigroups share an important property with groups, namely that there is a natural unary operation $x\mapsto x'$ which assigns to each element its unique inverse. Thus inverse semigroups are frequently viewed as algebras of type $\langle 2,1\rangle$ with the binary operation $\cdot$ being the semigroup multiplication and the unary operation ${}'$ being the natural inversion.

In the language of algebras of type $\langle 2,1\rangle$, a set of $n$ independent identities is an $n$-base for inverse semigroups, if those identities define the variety of inverse semigroups with the unary operation coinciding with the natural inversion.
The best known equational characterization of the variety of inverse semigroups is the following, due to B. Schein (\cite[Theorem 1.4]{Schein}):
\[
x(yz)=(xy)z
\qquad
(xy)'=y'x'
\qquad
x''=x
\qquad
xx'x=x
\qquad
(xx')(y'y)=(y'y)(xx')\,.
\]
This is not, strictly speaking, a base because
the identity $(xy)'=y'x'$ is dependent on the others, as Schein himself noted. Removing that identity turns out to give a $4$-base for inverse semigroups. One can also replace $x''=x$ with $x'xx' = x'$ thus giving another independent 4-base \cite[Theorem 3.2]{AM}. A $5$-base also due to Schein is given by keeping by the first four of the above identities and replacing $xx'y'y=y'yxx'$ with $xx'x'x=x'xxx'$.

So the next natural question was if there exists a $3$-base for inverse semigroups. Finding minimal axiom sets for classes of mathematical objects has, of course, a long and distinguished tradition in the area of automated deduction; see the bibliography in \cite{AM}. Not surprisingly, then, the affirmative solution to the $3$-base problem was found by the first-named author and Bill McCune using \textsc{Prover9} \cite{McCune}. That solution, which appeared in \cite{AM}, can be described as follows.

\begin{theorem}\label{final}
Consider the following identities in a binary operation and a unary operation:
\begin{align*}
(S_1) &\quad x(yz)=(xy)z & (S_4) &\quad x'xyy' = yy'x'x &\quad (S_7) &\quad (xy')z=x(z'y)' \\
(S_2) &\quad xx'x = x & (S_5) &\quad (x')' = x & (S_8) &\quad (xx')'x=x \\
(S_3) &\quad x'xx' = x' & (S_6) &\quad (xy)' = y'x' & (S_9) &\quad (xx')(yy')=(yy')(x'x)'
\end{align*}

As algebras $(S,\cdot,{}')$
of type $\langle 2,1\rangle$, where the unary operation coincides with the natural inversion, inverse semigroups can be defined  by each of the following sets of axioms:
\begin{enumerate}
   \item \emph{(}$S_1$\emph{)}, \emph{(}$S_2$\emph{)}, \emph{(}$S_4$\emph{)}--\emph{(}$S_6$\emph{)};
    \item \emph{(}$S_1$\emph{)}--\emph{(}$S_4$\emph{)};
    \item \emph{(}$S_1$\emph{)}, \emph{(}$S_2$\emph{)}, \emph{(}$S_4$\emph{)} and \emph{(}$S_5$\emph{)};
    \item \emph{(}$S_1$\emph{)}, \emph{(}$S_3$\emph{)}, \emph{(}$S_4$\emph{)} and \emph{(}$S_5$\emph{)};
    \item \emph{(}$S_7$\emph{)}--\emph{(}$S_9$\emph{)}.
  \end{enumerate}
\end{theorem}

The $3$-base ($S_7$)--($S_9$) was shortly followed by a much more elegant one \cite{AK} (again using automated deduction):
\[
\quad x(yz'')=(xy)z \qquad
\quad x=(xx')x \qquad
\quad (xx')(y'y) = (y'y)(xx') \,.
\]
Since then we have discovered even more elegant $3$-bases, including some which have associativity as one of the axioms. These will appear together with human proofs in \cite{AKP}.

\medskip

As was proved in \cite{AM}, there does not exist a $1$-base for inverse semigroups, that is, any $k$-base for inverse semigroups requires $k\geq 2$. So until now the main open problem in this line of research has been to prove or disprove the existence of a $2$-base.

\section{Early Attempts}

A naive approach to the problem of trying to find a base for some algebra is to generate many identities which can serve as candidates, and then to loop through those candidates, testing them against both a theorem prover and a finite model builder.

In order for the naive approach to have any chance of succeeding, one needs some \emph{a priori} idea of what properties a particular base must have. For instance, if there is a lower bound on the number of variables which are necessary or a lower bound on the length of the axioms (as measured by, say, symbol counting), then this can very helpful in reducing the size of the candidate pool.

Inverse semigroups form a uniform (or regular) variety of algebras, that is, in each identity, the set of variables appearing on the left of the identity equals the set of variables appearing on the right of the identity. As is well known, any consequence of uniform identities is also  uniform and hence any equational base must be uniform as well.   In addition, it is easy to see that one of the identities must have the form $u=x$ (a variable equal to a word), and so $u$ must depend only on $x$. Putting together these criteria with a few others, we have the following \cite{AKP}:

\begin{proposition}
A $2$-base for inverse semigroups must have the form $\{u(x) = x, s = t\}$ where $u(x)$ is a term involving only the variable $x$. The identity $s=t$ involves at least three variables, all of which occur in both $s$ and $t$. The unary operation occurs in $u(x)$ and in at least one of $s$ or $t$.
\end{proposition}

In the end, this description was not as helpful as one might hope.
In 2006, Bill McCune generated approximately half a million candidates for possible $2$-bases, but after looping through them as described above, none of them worked.

In 2010, we asked Bill if any of those sets of candidates could imply just
\[
\quad x(yz)=(xy)z \qquad
\quad x=x(x'x) \qquad
\quad (xx')(y'y) = (y'y)(xx')\,.
\]
A semigroup together with a unary operation ${}'$ satisfying just these identities will indeed be an inverse semigroup, but in general, the unary operation will not coincide with the natural inverse. As it turns out, the point was moot, because once again all candidates failed.

\medskip

We then tried another (unpublished) technique which has sometimes been helpful in finding smaller axiom sets for various structures. The idea is to take advantage of \textsc{Prover9}'s \emph{semantic guidance} feature. \textsc{Prover9} can take as input \emph{interpretations} which are models in the standard \textsc{Mace4} format. In the default given clause selection scheme for \textsc{Prover9} (assuming there are no hints), clauses are selected in the following order: $1$ oldest clause, $4$ lightest false clauses and $4$ lightest true clauses. Here ``false'' means false in all interpretations and ``true'' means true in at least one interpretation. The default interpretation (meaning that the interpretations list does not contain any models) is that false clauses are negative and true clauses are nonnegative.

The idea behind using semantic guidance to reduce axiom sets goes as follows. Suppose we have a theory with, say, $3$ independent formulas $A$, $B$ and $C$. We use \textsc{Mace4} and \textsc{Prover9} as follows to carry out the  following procedure and hence produce candidates that are more likely to be a 2-base:
\begin{enumerate}
\item \textsc{Mace4}: Find models in which $A$ and $B$ are true, but $C$ is false.
\item \textsc{Mace4}: Find models in which $A$ and $C$ are true, but $B$ is false.
\item \textsc{Prover9:} With all models in the interpretations list, generate consequences of $A$, $B$ and $C$ for a fixed amount of time (or fixed number of given clauses or any other criterion for stopping the job).
\item If \textsc{Prover9} generated any clauses marked false, they will be false in \emph{all} interpretations. Collect all such clauses.
\item For each false clause $D$, test if $\{A,D\}$ imply $B$ and $C$ using \textsc{Prover9} and \textsc{Mace4}.
\end{enumerate}

This procedure has been very useful to us in other situations, but it failed here. We used all the $3$-bases known to us (those mentioned above and others) and tried the procedure above with all of them, producing for each a $3G$ text file; but none of the candidates worked.

After this experiment and others like it, the first two authors were temporarily convinced that no $2$-base exists for inverse semigroups (as stated in the end of \cite{AK}).
This conviction was totally changed when we discovered $2$-bases for two proper varieties of inverse semigroups, namely commutative inverse semigroups and Clifford semigroups \cite{AKP}. This made the existence of a $2$-base for inverse semigroups much more plausible since it is unusual to have a non-$2$-based variety with $2$-based subvarieties.

\section{A $2$-base for inverse semigroups and its proofs}

As already noted, one identity in any $2$-base for inverse semigroups must have the form $u(x) = x$. Thus in our search for a $2$-base, we decided to fix this identity in the simplest way possible, namely $x(x'x) = x$ (or sometimes $(xx')x = x$). So what remained was to find a uniform identity which would imply associativity, the commutativity of idempotents and $x''=x$. After many tests using \textsc{Prover9}, we identified an identity that, together with $x(x'x)=x$, implies associativity and $x''=x$, but not commutativity of the idempotents: $(x  y)  z = ((y  z)'  x')'$.

Then we used a trick typical in the world of cancellative semigroups: we \emph{glued} together this identity with an identity expressing commutativity of idempotents. Somewhat surprisingly, one variant of this approach worked and yielded the following $2$-base:
\begin{align*}
& x = x(x'x) \\
& x  (x'  (y  (y'  ((z  u)'  w')'))) = y  (y'  (x  (x'  ((w  z)  u))))\,.
\end{align*}

What is remarkable is that making even small changes in the second identity, such as switching the roles of $y$ and $y'$, will fail to yield a $2$-base.

While the discovery of these identities required a mix of automated tools and human intuition, the verification that these identities form a $2$-base for inverse semigroups is straightforward for automated theorem provers. We initially checked it using \textsc{Prover9} and subsequently using \textsc{Waldmeister} \cite{Hillenbrand} and \textsc{E} \cite{Schulz}.

The first proof found by any of these tools is quite lengthy and complex. We have been able to generate simpler proofs using \textsc{Prover9}'s hints feature \cite{Veroff}. The basic idea is to find an initial proof, use it to generate hints, run \textsc{Prover9} again to find a new proof and repeat. We do not have formal criteria for simplicity (\emph{cf.} \cite{TW}), but rather look for proofs which find a balance between length, depth, the number of variables used and so on.

The underlying idea here follows a dictum oft repeated by McCune: ``It's all about the given clause." Changing which clauses will be selected as given in subsequent runs changes the search space, thus leading to the possibility that new, simpler proofs will be found.

For example, in practice, when \textsc{Prover9} is fed hints from its first successful proof, it rarely reproduces that same proof in the second run. Typically, the proofs of key steps are found more quickly so that the whole second proof is shorter and perhaps simpler than the first proof.

At the same time, in order to reduce the complexity of the proof further, we gradually reduce the maximum weight parameter to try to find a proof using shorter and shorter clauses. (Note that by default, hint matchers are exempt from \textsc{Prover9}'s various limits on kept clauses, so this particular strategy requires setting a flag.) We reduce the maximum clause weight until just before a proof can no longer be found. Of course, there is no reason to suspect that a proof which has the minimum possible maximum clause weight will necessarily be the simplest proof, so slight increases in the maximum weight from its minimum sometimes yield better results.

There is a law of diminishing returns in this in that eventually a point is reached where the effort involved in further simplifying a proof goes beyond any additional minor simplifications that might occur. In our case, we eventually found a 92 step proof of $x'' = x$, a 29 step proof of associativity (assuming $x''=x$) and a 26 step proof of $xx'y'y=y'yxx'$ (assuming $x''=x$ and associativity). Here we are referring to \textsc{Prover9}'s own reporting of proof length which, in proofs with demodulations, counts only the primary inferences. If we count rewrites as inferences, the proof lengths are 274, 82 and 93, respectively.  However, we find that proofs treating rewrites as secondary inferences are often easier to follow for a human reader than those which treat primary paramodulations and rewrites on an equal footing.

Of course the three proofs, with or without demodulations, share many inferences in common. Interestingly, adjoining the three goals together into a single goal and getting a proof in one run leads to a more complicated proof. This is probably because when separated, the second proof takes early advantage of the additional assumption $x''=x$, and the third proof takes advantage of both that and associativity.

In a subsequent publication \cite{AKP}, we will present humanized proofs that our pair of identities forms a $2$-base for inverse semigroups, as well as some interesting $3$-bases and some $2$-bases for certain subvarieties of the variety of inverse semigroup. For now, we present the three proofs, lightly edited. The justification steps list only the parent clauses, not how they are specifically used to make the inference.

{\tiny{
\begin{verbatim}
% Proof 1 at 0.23 (+ 0.01) seconds.
% Length of proof is 92.
% Level of proof is 36.
% Maximum clause weight is 49.000.
% Given clauses 70.

1 x'' = x # label(goal).
2 x * (x' * x) = x.
3 x * (x' * (y * (y' * ((z * u)' * w')'))) = y * (y' * (x * (x' * ((w * z) * u)))).
4 c1'' != c1.  [1].
5 ((x * y)' * z')' * (((x * y)' * z')'' * (u * (u' * ((z * x) * y)))) = u * (u' * ((x * y)' * z')').  [2,3].
6 x * (x' * (y * (y' * (((z' * z) * u)' * z')'))) = y * (y' * (x * (x' * (z * u)))).  [2,3].
7 ((x * y) * z) * (((x * y) * z)' * (u * (u' * ((y * z)' * x')'))) = u * (u' * ((x * y) * z)).  [2,3].
8 ((x * y)' * z')' * (((x * y)' * z')'' * ((z * x) * y)) = ((z * x) * y) * (((z * x) * y)' * ((x * y)' * z')').  [2,7].
9 ((x * y)' * z')' * (((x * y)' * z')'' * (((z * x) * y) * (((z * x) * y)' * ((x * y)' * z')'))) = ((x * y)' * z')'.  [2,5,8].
10 ((x * y)' * z')' = (z * x) * y.  [9,3,8,7,2].
11 x * (x' * (y * (y' * (z * u)))) = y * (y' * (x * (x' * (z * u)))).  [6,10,2].
12 (x * y) * (y' * y) = (y' * x')'.  [2,10].
13 (((x * y) * z) * u')' = (u * (y * z)') * x'.  [10,10].
14 (x * (y * z)) * (z' * z) = ((z' * y')'' * x')'.  [12,10].
15 (x * ((y' * y) * z)') * y' = ((y * z) * x')'.  [2,13].
16 (x * (y' * y'')'') * y' = (y * x')'.  [2,15,12].
17 (x * y) * ((y' * y)' * (y' * y)) = (x * (y * y')) * y.  [2,14,10,10].
18 ((x' * x'')'' * x')' = x.  [2,14,2].
19 ((x' * x'')'' * (y * x)')' = (x' * y')' * (x' * x).  [12,14].
20 (((x' * x'')'' * x') * y) * z = ((y * z)' * x)'.  [18,10].
21 ((x * y) * ((x * y)' * x)) * y = x * y.  [10,18,10].
22 ((x * x')'' * x)' = (x' * x'')'' * x'.  [18,18,18,18].
23 (x'' * x''')'' * x'' = x.  [18,22].
24 (x * ((y * z) * ((y * z)' * y))) * z = (x * y) * z.  [21,10,10].
25 (((x * y)' * z') * (((z * x) * y) * (x * y)')) * z' = (x * y)' * z'.  [10,21].
26 (x * (y * ((z * u) * ((z * u)' * z)))) * u = (x * (y * z)) * u.  [24,10,10].
27 (x' * (x * x')) * x = x' * x.  [17,2].
28 (x * (y' * (y * y'))) * y = (x * y') * y.  [27,10,10].
29 (x' * (x' * (x * x'))')' = (x' * x'')'.  [27,12,12].
30 (x' * (y * (x' * (x * x')))')' = (x' * (y * x')')'.  [28,12,12].
31 (x * y') * (y' * (y * y'))' = (x * y') * y''.  [29,10,10].
32 x' * (x' * (x * x'))' = x' * x''.  [29,23,29,29,23].
33 (x * (y * z')) * (z' * (z * z'))' = (x * (y * z')) * z''.  [31,10,10].
34 x' * (y * (x' * (x * x')))' = x' * (y * x')'.  [30,23,30,30,23].
35 ((x * (y' * (y * y'))) * y'')' = ((x * y') * y'')'.  [34,15,15].
36 (x * (y' * (y * y'))) * y'' = (x * y') * y''.  [35,23,35,35,23].
37 ((x * (y * z)) * ((y * z)' * y)) * z = ((x * y) * z) * ((y * z)' * (y * z)).  [10,19,10,10].
38 (x * (y * (z * u'))) * (u' * (u * u'))' = (x * (y * (z * u'))) * u''.  [33,10,10].
39 x * (x' * (y * (y' * z))) = y * (y' * (x * (x' * z))).  [2,11,2].
40 x * (x' * (y * (y' * x))) = y * (y' * x).  [2,39].
41 ((x * y) * ((x * y)' * (z * (z' * x)))) * y = (z * (z' * x)) * y.  [39,26].
42 (x'' * x') * ((x'' * x')' * (x * x')) = x * x'.  [2,40,2].
43 (x * y) * (y' * (z * (z' * y))) = (x * z) * (z' * y).  [40,10,10].
44 (x * (x' * y)) * (y' * x) = y * (y' * x).  [40,14,2,10,10].
45 ((x * (x' * y)) * y') * x = (y * y') * x.  [44,12,12,10,10].
46 (x' * ((x * (x' * y)) * y')')' = (x' * (y * y')')'.  [45,12,12].
47 ((x * y) * z) * ((y * z)' * (y * z)) = (y' * x')' * z.  [43,26,37,12].
48 ((x * y') * (y * y')) * y = (y'' * x')' * y.  [47,17].
49 ((x * y)' * z)' * ((x * y)' * (x * y)) = (x' * z)' * y.  [20,47,22,23].
50 ((x * y') * y'') * (((y' * (y * y')) * y'')' * ((y' * (y * y')) * y'')) = ((x * y') * (y * y')) * y''.  [36,47,10].
51 (x' * ((y * x') * (x * x'))')' = (x' * (x'' * y')'')'.  [48,12,12].
52 x' * ((x * (x' * y)) * y')' = x' * (y * y')'.  [46,23,46,46,23].
53 (x'' * x''')'' * (((x'' * x''')' * x) * x''')' = (x'' * x''')'' * (x'' * x''')'.  [23,52].
54 x * (((x * y) * y')' * x)'' = x * (y * y')'.  [20,52,22,23,22,23,22,23].
55 x * ((y' * (y * y')) * y'')' = x * (y' * y'')'.  [33,54,38,36,54].
56 (x'' * y)' * (x' * (x * x'))' = (x'' * y)' * x''.  [32,49,32,32,49].
57 ((x' * x'') * (((x' * (x * x')) * x'')' * x')) * x'' = (x' * (x * x')) * x''.  [55,22,10,36,10,23].
58 x''' * (x' * (x * x'))' = x''' * x''.  [2,56,2].
59 (x' * (x * x'))' * ((x' * (x * x'))'' * x'') = x''.  [58,40,2,58,2].
60 (((x' * (x * x'))' * y) * (((x' * (x * x'))' * y)' * x'')) * y = x'' * y.  [58,41,2,58,2].
61 x' * ((y * x') * (x * x'))' = x' * (x'' * y')''.  [51,23,51,51,23].
62 ((x' * (x * x')) * x'')' = (x' * x'')'.  [61,15,16].
63 (x' * (x * x')) * x'' = x' * x''.  [57,62,21].
64 ((x * y') * (y * y')) * y'' = (y'' * x')' * y''.  [50,63,63,47].
65 ((x' * (x * x'))' * y)' * x'' = (x'' * y)' * x''.  [63,49,63,63,49].
66 (((x' * (x * x'))' * y) * ((x'' * y)' * x'')) * y = x'' * y.  [60,65].
67 (x' * (x * x'))' * y' = x'' * y'.  [33,25,64,66].
68 (x * y') * (y * y') = (y'' * x')'.  [67,10].
69 ((x' * (x * x'))' * (x' * x''))' * y' = (x' * (x * x'))'' * y'.  [33,67,63].
70 (x' * (x * x'))'' * x'' = x''' * x''.  [67,65,33,63,67].
71 (x' * (x * x'))' * (x''' * x'') = x''.  [59,70].
72 ((x' * (x * x'))' * y)' = (x'' * y)'.  [68,20,22,23].
73 (x' * (x * x'))'' * y' = x''' * y'.  [69,72,67].
74 (x' * (x * x'))' = x''.  [73,2,58,71].
75 x' * (x * x') = x'.  [73,23,74,74,23].
76 (x'' * y)' * (x'' * x') = (x'' * y)' * (x * x').  [75,49,75,75].
77 x'' * x' = x * x'.  [76,2,42].
78 (x''' * y)' * x' = (x' * y)' * x'.  [76,49,77,77,49].
79 (x'' * x''')'' * (((x'' * x''')' * x) * x''')' = (x'' * x''') * (x'' * x''')'.  [53,77].
80 (x'' * x''')' = (x * x')'.  [77,12,77,68,77].
81 ((x * x') * (x * x')')' = (x * x')' * (x * x').  [77,19,80,77,80,77].
82 (x * x')'' * (((x * x')' * x) * x''')' = (x'' * x''') * (x * x')'.  [79,80,80,80].
83 (x * x')'' * x'' = x.  [23,80].
84 (x'' * x''') * (x * x')' = (x * x') * (x * x')'.  [83,52,82,80,77].
85 ((x * x')' * (x * x'))' * (x * x')'' = x'' * x'''.  [80,83,84,81,80].
86 (x' * (x'' * x'''))' * x' = x'''' * x'.  [2,78,77].
87 x'' * x''' = x * x'.  [81,83,85].
88 x'''' * x' = x * x'.  [86,87,75,77].
89 x''' * (x * x') = x'.  [87,40,75,88,87,75].
90 x''' = x'.  [87,75,89].
91 x'' = x.  [87,83,90,83].
92 $F.  [91,4].
============================== end of proof ==========================

% Proof 1 at 0.04 (+ 0.00) seconds.
% Length of proof is 29.
% Level of proof is 14.
% Maximum clause weight is 47.000.
% Given clauses 25.

1 (x * y) * z = x * (y * z) # label(goal).
2 x * (x' * x) = x.
3 x * (x' * (y * (y' * ((z * u)' * w')'))) = y * (y' * (x * (x' * ((w * z) * u)))).
4 x'' = x.
5 (c1 * c2) * c3 != c1 * (c2 * c3).  [1].
6 x * (x' * (y * (y' * ((z * u) * (u' * u))))) = y * (y' * (x * (x' * (u' * z')'))).  [2,3].
7 ((x * y)' * z')' * (((x * y)' * z') * (u * (u' * ((z * x) * y)))) = u * (u' * ((x * y)' * z')').  [2,3,4].
8 ((x * y) * z) * (((x * y) * z)' * (u * (u' * ((y * z)' * x')'))) = u * (u' * ((x * y) * z)).  [2,3].
9 x' * (x * x') = x'.  [4,2].
10 ((x * y) * (y' * y)) * (((x * y) * (y' * y))' * (z * (z' * (y' * x')'))) = z * (z' * ((x * y) * (y' * y))).  [2,6].
11 x * (x' * (y' * (y * ((z * u) * (u' * u))))) = y' * (y * (x * (x' * (u' * z')'))).  [4,6,4].
12 x' * (x * (y' * (y * ((z * u) * (u' * u))))) = y' * (y * (x' * (x * (u' * z')'))).  [4,11,4].
13 ((x * y)' * z')' * (((x * y)' * z') * ((z * x) * y)) = ((z * x) * y) * (((z * x) * y)' * ((x * y)' * z')').  [2,8,4].
14 ((x * y)' * z')' * (((x * y)' * z') * (((z * x) * y) * (((z * x) * y)' * ((x * y)' * z')'))) = ((x * y)' * z')'.  [2,7,4,13].
15 ((x * y) * (y' * y)) * (((x * y) * (y' * y))' * (y' * x')') = (y' * x')' * ((y' * x') * ((x * y) * (y' * y))).  [2,10,4].
16 (x * y) * (y' * y) = (y' * x')'.  [2,10,15,11,15,12,9,9].
17 (x' * y)' * ((x' * y) * (z * (z' * (x' * y)'))) = z * (z' * (x' * y)').  [4,10,16,4,16,4,4,16,4].
18 ((x * y) * z) * (((x * y) * z)' * ((y * z)' * x')') = ((y * z)' * x')'.  [14,17].
19 (x * y') * (y * y') = (y * x')'.  [4,16,4].
20 ((x * y)' * z')' = (z * x) * y.  [18,8,18,2].
21 ((x * y) * z)' = (y * z)' * x'.  [20,4].
22 ((x * y)' * z)' = (z' * x) * y.  [4,20].
23 ((x * y)' * x) * x' = (x * y)'.  [2,21,21,4].
24 (x' * (y * z)) * u = (y' * x)' * (z * u).  [21,22,21,4].
25 (x * (y * z))' = (y * z)' * x'.  [23,22,4,21,24,4,21,19,4].
26 (x * y)' = y' * x'.  [2,25,2].
27 (x' * y) * z = x' * (y * z).  [22,26,26,26,4,4].
28 (x * y) * z = x * (y * z).  [4,27,4].
29 $F.  [28,5].
============================== end of proof ==========================

% Proof 1 at 0.02 (+ 0.00) seconds.
% Length of proof is 26.
% Level of proof is 11.
% Maximum clause weight is 34.000.
% Given clauses 21.

1 (x * x') * (y' * y) = (y' * y) * (x * x') # label(goal).
2 x * (x' * x) = x.
3 x * (x' * (y * (y' * ((z * u)' * w')'))) = y * (y' * (x * (x' * ((w * z) * u)))).
4 x'' = x.
5 (x * y) * z = x * (y * z).
6 (c2' * c2) * (c1 * c1') != (c1 * c1') * (c2' * c2).  [1].
7 c2' * (c2 * (c1 * c1')) != c1 * (c1' * (c2' * c2)).  [6,5,5].
8 x * (x' * (y * (y' * ((z * u)' * w')'))) = y * (y' * (x * (x' * (w * (z * u))))).  [3,5].
9 x' * (x * x') = x'.  [4,2].
10 x * (x' * (x * y)) = x * y.  [2,5,5].
11 x * (x' * ((y * (z * u))' * x')') = x * ((z * u)' * y')'.  [8,8,4,4,4,4,10,4,4,4,4,10,10].
12 x' * (x * (x' * y)) = x' * y.  [9,5,5].
13 x * (x' * ((y * z)' * u')') = x * (x' * (u * (y * z))).  [10,8,12].
14 x * (x' * (y * (y' * (z * u)))) = y * (y' * (x * (x' * (z * u)))).  [10,8,13,10].
15 x * ((y * z)' * u')' = x * (u * (y * z)).  [11,13,10].
16 x * (y * z')' = x * (z * y').  [9,15,4,9].
17 x * (y * z)' = x * (z' * y').  [4,16].
18 x * ((y * z)' * u) = x * (z' * (y' * u)).  [17,5,5,5].
19 x * (y * (y' * (x' * (z * (z' * (x * y)))))) = z * (z' * (x * y)).  [2,14,18,5].
20 x * (x' * (y * (y' * (x * z)))) = y * (y' * (x * z)).  [10,14].
21 x * (x' * (y' * (y * (x * z)))) = y' * (y * (x * z)).  [4,20,4].
22 x' * (x * (y * (y' * x'))) = y * (y' * x').  [9,20,4,9].
23 x * (y * (z * (z' * (y' * y)))) = x * (y * (z * z')).  [22,17,17,18,4,4,4,18,18,18,4,4].
24 x * (y * (y' * (x' * x))) = x * (y * y').  [23,2,18,18,18,4,21,18,4,10,5,5,5,5,10,21,10].
25 x' * (x * (y * y')) = y * (y' * (x' * x)).  [19,12,4,24].
26 $F.  [25,7].
============================== end of proof ==========================
\end{verbatim}}}

\section{Problems}

We conclude with a couple of open problems which are suitable for
investigation by means of automated deduction.

\begin{prob}
Does there exist a $2$-base with fewer than $5$ variables?
\end{prob}

\begin{prob}
Taking the first identity $x(x'x)=x$ as fixed, is there a $2$-base with a
shorter second identity?
\end{prob}

We have found $2$-bases for certain subvarieties of inverse semigroups, such as
commutative inverse semigroups and Clifford semigroups; these will appear in \cite{AKP}.
But there are other subvarieties, such as strict inverse semigroups \cite{Petrich} for
which we have not found $2$-bases.

\begin{prob}
Find $2$-bases for other varieties of inverse semigroups.
\end{prob}

\begin{acknowledgment}
The first author was partially supported by FCT and FEDER, Project POCTI-ISFL-1-143 of Centro de Algebra da Universidade de Lisboa, and by FCT and PIDDAC through the project PTDC/MAT/69514/2006.

The third author was supported by a University of Manitoba research leave grant during 2012-13
and he thanks the University of Manitoba for sanctioning Research-Study Leave.
\end{acknowledgment}

\end{document}